\documentclass[a4paper,12pt]{article}
\usepackage{hyperref}
\usepackage{amsthm, amsmath, amsfonts, amssymb}
\usepackage{enumerate}
\usepackage[initials, alphabetic]{amsrefs}

\newtheorem{thm}{Theorem}[section]

\numberwithin{equation}{section}

\theoremstyle{definition}
\newtheorem{rem}{Remark}[section]

\theoremstyle{remark}

\renewcommand{\P}{{\mathbf P}}
\newcommand{\R}{\mathbf{R}}
\newcommand{\cN}{\mathcal{N}}
\newcommand{\B}{\mathcal{B}}
\newcommand{\E}{\mbox{\boldmath $E$}}
\newcommand{\bW}{\mathbf W}
\newcommand{\bxi}{\xi}
\newcommand{\F}{\mathcal{F}}
\newcommand{\Fo}{\F}

\begin{document}
\title{
What is the natural scale for a L\'evy process 
in modelling term structure of 
interest rates?
}

\author{Jir\^o Akahori\thanks{This research is 
partially supported by Open Research Center Project for Private Universities: matching fund subsidy from MEXT, 
2004-2008
and also by Grants-in-Aids for Scientific Research (No. 18540146)
from the Japan Society for Promotion of Sciences.} 
\,and Takahiro Tsuchiya \\
Graduate School of Mathematics, Ritsumeikan University \\
1-1-1 Nojihigashi, Kusatsu, Shiga, 525-8577, Japan \\
E-mail: akahori@se.ritsumei.ac.jp, suci@probab.com}

\date{}
\maketitle

\begin{abstract}                       
This paper gives 
examples of explicit arbitrage-free 
term structure models with L\'evy jumps 
via state price density approach.
By generalizing quadratic Gaussian models,
it is found that the probability density function
of a L\'evy process is a "natural" scale 
for the process to be the state variable of a market.
\end{abstract}

\thanks{\noindent {\bf Keywords}. State price density approach, term structure models, 
Shirakawa model, L\'evy process, Probability density.}

\thanks{\noindent {\bf 2000 Mathematics Subject Classification}:
91B70, 60G52}
\noindent\thanks{{\bf Journal of Economic Literature 
Classification System}:G12}

\section{Introduction}
\subsection{Literature review}
In the classical Black-Scholes economy, 
the exponential of a Brownian motion (with drift), i.e.,
\begin{equation*}
S_t = S_0 \exp (a B_t + b t), 
\end{equation*}
is used to model stock prices. 
When one takes jumps into account, 
extending the Black-Scholes scale 
naturally leads to the modelling economic factors by exponential 
of a L\'evy process, say:
\begin{equation*}
X_t = X_0 \exp  Z_t,
\end{equation*}
where $ Z $ is a L\'evy process. 
Such a model is often called 
"exponential-L\'evy" type, 
and has been widely used in financial modelling
since Merton's jump diffusion model \cite{Mer} appeared.
The variance gamma model by Madan and his co-authors 
(see e.g. \cite{MCC}), 
and the hyperbolic model by Eberlein and his authors 
(see e.g. \cite{EbPr}),
are two major representatives, 
but there are still many others including \cite{B-NO} and 
\cite{B-NL}. For further references, 
see \cite{Miy} or \cite[Chapter 8--11]{CoTa}. 

In the context of interest rate models, 
the exponential scale is also common
(though the scale is in general
not time-homogeneous any more);
e.g., bond prices of Vasicek's model \cite{Vas} 
or Gaussian HJM \cite{HJM}, 
LIBOR rates of the BGM's LIBOR market model \cite{BGM}, etc.  
There are also exponential-L\'evy type interest rate models with jumps;
e.g. Shirakawa's extension of Gaussian HJM \cite{Shi}, 
Chiarella-Nikitopoulos \cite{Ch-Ni}, 
Eberlein-Raible's L\'evy term structure models \cite{EbRa}, 
or Albeverio-Lytvynov-Mahnig \cite{AlLyMa}, etc, 
etc (see also section \ref{LLTSM} in the present paper). 

There are two major substitutes in the continuous-path cases; 
the affine class and the quadratic class.
Roughly speaking, we say the model is ``(exponential-)quadratic" 
(w.r.t. a state process $ Z $) 
if the price at time $ t $ of the bond 
with maturity $ T $ is given by 
\begin{equation*}
P^T_t = \exp \left\{ \langle A_2(t,T) Z_t, Z_t \rangle 
+ \langle A_1 (t,T), Z_t \rangle + A_0 (t,T) \right\},
\end{equation*}
where 
$ A_2 (t,T) \in \R^n \otimes \R^n $ is a 
non-zero symmetric matrix, 
$ A_1 (t,T) $ is in $ \R^n $, and $ A_0 (t,T) \in \R $. 
If $ A_2 \equiv 0 $ then the model is called 
``(exponential-)affine". 

Affine class w.r.t. a L\'evy process is nothing but
the above-mentioned exponential-L\'evy class. 
When we refer to affine class, 
the state process is usually a Markov process 
whose state space is a convex cone in $ \R^n $ 
(and hence symmetric processes are excluded).
The class contains 
Cox-Ingersoll-Ross model \cite{CIR} 
and Duffie-Kan's multi-dimensional generalization \cite{DuKa}.  
The jump-diffusion extension of the affine class 
was fully done by Duffie-Filipovi\'c-Schachermayer \cite{DFS} 
by generalizing Kawazu-Watanabe theorem \cite{KaWa}.
It is also shown by 
Filipovi\'c-Teichmann \cite{Fi-Te} that 
the "finite-dimensional realization" is limited to the affine class. 

The quadratic class, which can be 
embedded into the (degenerate) affine class, is
gaining much popularity among economists and practitioners
mainly due to its tractability,
(see e.g. \cite{ARG})
or partially explained from 
a mathematical background; it has a 
rich mathematical structure (see \cite{AkaHar}). 
Contrary to the affine case, however, no jumps are allowed 
in the quadratic class 
as Chen-Filipovi{\'c}-Poor showed in \cite{CFP} 
within Markovian framework (though
Levendorski{\u\i} \cite{Lev} 
gives a jump-type extension of quadratic models
``in the direction of affine class").

The main aim of the present paper 
is to give a description of 
a L\'evy extension of the quadratic models 
in a totally different way from \cite{CFP} and \cite{Lev}.

\subsection{Our result}
The above story of the quadratic term structure models
tells us that exponential-quadratic functions
are "natural" scales (at least in modelling interest rates)
for Brownian motions 
but not for pure-jump L\'evy processes. 
Then what scale is natural for a L\'evy process? 

To answer this question we start from the observation 
that within the exponential-quadratic scale
lies the probability density of the Gaussian distribution. 
This observation naturally leads to asking if 
its probability density is a {\em natural} scale for 
a general L\'evy process or not. 
To be more precise, the question is: 
given a L\'evy process $ Z $, the market consists of
\begin{equation*}
P^T_t = p(A(t,T), Z_t) \quad (\text{or something like this})
\end{equation*}
where $ p (t, x) = \P (Z_t \in dx | Z_0 = 0 )/dx $
can be consistent? 
or in particular, arbitrage-free?

An answer to this question is given as Theorem \ref{mainT}.
It says that if the instantaneous forward rate 
$ f (t,T) := - \partial_T \log P^T_t $ is given by 
\begin{equation*}
-\partial_T \log p( \lambda_T + T-t , Z_t ) 
\end{equation*}
for some continuous $ \lambda :[0,\infty) \to [0,\infty) $,
then the market is arbitrage-free (Theorem \ref{mainT}). 

To construct arbitrage-free interest rate model, 
we rely on so-called {\em state price density}, 
or {\em pricing kernel} approach \cite[pp371--]{BriMer}, 
which was initiated by Constantinides \cite{Con}
and later developed by Rogers \cite{Rog} and 
Hughston and his co-authors \cite{FleHug, BroHug, HugRam} 
(see also the textbook by 
Hunt and Kennedy \cite{HunKen}).

\subsection{Organization of the present paper}
We will start from a brief survey of the approach 
(section \ref{review}), and then give
two important classes; Gaussian (section \ref{GaussHJM}) 
and quadratic Gaussian (section \ref{QTSM}).
A direct jump-type extension which we call
{\em generalized Shirakawa model} 
will be given in section \ref{LLTSM}.
In section \ref{aSTSM1}, 
Theorem \ref{mainT} and its proof will be given. 
Examples based on our new framework will be presented 
in section \ref{exs}.

\section{The state price density approach 
to interest rate modelling}

\subsection{Review}\label{review}
In principle,  
a strictly positive process
$ \{ \pi_t\} $
is a {\em state price density} with respect to a market on a  
filtered probability space $ (\Omega, \F,  \P ,\{ \F_t \}) $
if for any asset indexed by $ i \in I $ that generates 
$ \{ D^i_s\} $ cash flows in the future, 
its market price at time $ t $ is given by 
\begin{equation}\label{CFF}
S^i_t = \pi_t^{-1} 
\E^\P [ \int_{t+}^\infty \pi_s dD^i_s | {\cal F}_{t} ], 
\end{equation}
or for any $ T  >  t $, 
\begin{equation*}
S^i_t = \pi_t^{-1} \E^\P [ \pi_T S^i_T + \int_{t+}^{T-} \pi_s dD^i_s | {\F}_{t} ].
\end{equation*}
In other words, 
$ \pi_T/\pi_t $ ($\times$ probability density with respect to $ \P $)
gives (random) discount factor of a (random) cash flow at time $ T $.
In particular, if we denote by $ P^{T}_{t} $ 
the market value at time $ t $
of zero-coupon bond with maturity $ T $, 
we have 
\begin{equation}\label{SPD}
P^{T}_{t}= \pi_t^{-1} \E^\P [ \pi_T| {\F}_{t} ].
\end{equation}

In the famous text book by Duffie \cite{Duf},
an arbitrage-free market is characterized by 
the existence of a  state price density.
This {\em duality}, which comes from  
the Hahn-Banach theorem and its variants,
is proven for fully discrete markets (finite dimensional cases) and 
for Brownian markets (Brownian filtration cases). 
In more general cases one can extend the duality though
she/he needs to be careful about 
the exact meaning of the {\em arbitrage}.
It depends on what kind of assets are traded 
and what kind of trading strategies are admissible in the market.
To determine how far we can extend is, however, 
out of the scope of the present paper\footnote{Extensive studies 
in this direction are found 
in \cite{DS95, DS99, DSText}, and 
for the cases of bond markets
\cite{BMKR}, \cite{BjKaRu}, and \cite{DeD}
to name a few. }. 
Here we just assume that the space of the value process of 
trading strategies is orthogonal 
to $ \pi $ in the $ L^2 $ space of stochastic processes, 
meaning that we presume that the existence of 
a state price density implies the market is arbitrage-free
(though we leave the problem of market completeness behind). 
This assumption is robust because at least it is fulfilled if we admit 
{\em simple strategies}\footnote{by which we mean 
those strategies which remain constant for a short time interval.} only. 

From a perspective of modeling term structure of interest rates,
the formula (\ref{SPD}) says that, given a filtration, 
each strictly positive process $ \pi $ generates 
an {\em arbitrage-free} interest rate model. 
On the basis of this observation, we will construct arbitrage-free
interest rate model.

\subsection{Examples}
In practice, an explicit formula for the bond price (and hence
interest rates) is desirable;
it gives a parameterization of 
the entire term structure, and it becomes applicable
to, say, duration-based hedging. 
The more we have explicit formulas for present values
$ S^i $ for the cash flow $ D^i $ through (\ref{CFF})
the better the model becomes.

Below we give two classical {\bf T}erm {\bf S}tructure 
{\bf M}odels (TSMs for short) which exhibit no jumps.
\
\subsubsection{Gaussian TSMs}\label{Gaussian}
Let $ \bW_t $ be a standard Brownian motion taking values in $\R^d$
and $ h (t,s) $ be an element of $ H \otimes H $ where 
\begin{equation*}
H= \{ h: [0,\infty) \to \R^d, \mbox{absolutely continuous s.t.}\,
\dot{h} \in L_{\mathrm{loc}}^2 [0,\infty) \}.
\end{equation*}
Given an initial data $ T \mapsto P_0^T $, define 
\begin{equation}\label{GaussSPD}
 \pi_t= P^t_0 
 \exp \left\{ \int_0^t \langle h_s (t,s) ,\,d\bW_s \rangle_{\R^d}
 - \frac{1}{2} \int_0^t {| h_s (t,s) |}^2 \,ds  \right\}. 
\end{equation}
Then, by an easy manipulation we have
\begin{equation}\label{GaussB}
\begin{split}
& P^T_t = \frac{P^T_0}{P^t_0}
 \exp \bigg\{ \int_0^t \langle h_{s} (T,s) 
 - h_{s} (t,s) ,\,d\bW_s \rangle_{\R^d} \\
& \hspace{3cm}
- \frac{1}{2} \int_0^t ({| h_s (T,s) |}^2 
- {| h_s (t,s) |}^2) \,ds  \bigg\},
\end{split}
\end{equation}
where $ h_s $ stands for the partial derivative with respect to 
the latter variable. 

Note that this class covers 
so-called Gaussian Heath-Jarrow-Morton models
(see e.g. \cite[319--]{MusRut}). 
In fact, we have
\begin{equation}\label{GaussHJM}
\begin{split}
&f(t,T) := - \partial_T \log P^T_t \\
&=  f(0,T) 
+ \int_0^t \langle - h_{T,s} (T,s)  ,\,d\bW_s \rangle_{\R^d}
+ \int_0^t \langle h_{T,s} (T,s) , h_s (T,s) \rangle_{\R^d} \,ds \\
&= f(0,T) 
+ \int_0^t \langle - h_{T,s} (T,s)  , d\bW_s 
- h_s (T,s) \,ds \rangle_{\R^d}.
\end{split}
\end{equation}
The last expression shows that 
$ \{ f (t,T) \}_{t \leq T} $ is a martingale 
under the so-called {\em forward measure} $ \P^T $ defined by
\begin{equation*}
\begin{split}
d\P^T /d\P 
&= \frac{P^T_t }{\E [P^T_t]}   \\
&= \exp \left\{ -  \int_0^T \langle h_{s} (T,s) 
,  d\bW_s \rangle_{\R^d} 
- \frac{1}{2} \int_0^t |h_{s} (T,s) |^2\,ds  \right\}. 
\end{split}
\end{equation*}

\subsubsection{Quadratic Gaussian TSMs}\label{QTSM}
Let $ A $ be a continuous map on $ \R_+ $
taking values in 
the set of all positive definite $ d \times d $-symmetric matrices,
and $ k : \R_+ \to \R $ be a continuous map.
Define
\begin{equation*}
\pi_t (x) = \exp \left\{ 
-  \langle A_t x, x \rangle_{\R^d} + k_t  
\right\}, \quad (x \in \R^d). 
\end{equation*}
Then we have, for a $ d $-dimensional Wiener process $ \bW$
starting from the origin, 
\begin{equation}\label{QTSM1}
\begin{split}
&P^T_t = \E [\pi_T (\bW_T) | \sigma (\bW_s; s \leq t ) ]/ 
\pi_t (\bW_t) \\
&=  \{ \det (2 (T-t) A_T + I )\}^{-1/2} \cdot\\
& \cdot 
\exp \left( -\langle  \left\{ A_T - A_t - 2(T-t) ( 2(T-t) A_T + I )^{-1} \right\}
\bW_t, \bW_t \rangle + (k_T-k_t)  \right), 
\end{split}
\end{equation}
where $ I $ is the unit matrix. 

For the derivation of (\ref{QTSM1}), see Appendix.


\begin{rem}
There are several ways to introduce QTSMs. 
The above is just an illustration. 
For details see e.g. \cite{CFP} and \cite{AkaHar}. 
\end{rem}

\subsubsection{Generalized Shirakawa TSMs}\label{LLTSM}

Let $ p $ be a stationary Poisson point process
on a measurable space $ (E, \B_E) $.
Then its counting measure defined by 
\begin{equation*}
\cN_p ((s,t], A) := \sharp \{ u \in (s,t]; p(u) \in A \}; \,\,
(s<t,\, A \in \B_E) 
\end{equation*}
is a stationary Poisson random measure; i.e. 
for mutually disjoint $ B_1,...,B_n \in \B_{\R_+ \times E} $, 
the random variables $ \cN_p (B_j) $'s are mutually independent
and Poisson distributed, and for $ A \in \B_E $ and $ s< t $,  
\begin{equation*}
\P (  \cN_p( (s,t], A) =n ) = \{ (t - s) \nu(A)\}^n 
\frac{e^{- (t-s) \nu(A)}}{n!},
\end{equation*}
where $ \nu $ is a $ \sigma $-finite measure  on $ (E, \B_E) $. 
For a detailed instruction, see e.g. \cite{IkeWat}. 

Let $\delta$ be a Borel function from 
$\R_+ \times \R_+ \times E $ to ${\bf R}$ 
such that (the equivalence class of)
$ \delta (t, \cdot, \cdot) $ 
is in $ \cap_{p \geq 1} L^p ( dt \otimes \nu ) $ and
$ \sup_{p \geq 1} \Vert \delta \Vert^p_{L^p}  < \infty $. 

For this $ \delta $ and $ h \in H \otimes H $, define
\begin{equation*}
Z_t = \int_0^t h_s (t,s) d\bW_s 
- \frac{1}{2} \int_0^t |h_s (t,s)|^2 \,ds
+ \int_0^t \int_E \delta (t,s,x) \,\cN_p (ds, dx),
\end{equation*}
where $ H $, $ h_s $, and $ \bW $ are as in Example \ref{Gaussian},
and $ \bW $ is independent of $ p $ (or equivalently, $ \cN_p $). 

By the assumption on $ h $ and $ \delta $, 
we have
\begin{equation*}
\begin{split}
\E [e^{Z_t}] &= \exp \left\{
\int_0^t \int_E (e^{ \delta (t,s,x)}- 1) \nu (dx) ds \right\} \\
&= \exp \left\{ \sum_{p=1}^\infty 
\int_0^t \int_E \{ \delta (t,s,x)\}^p \nu (dx) ds \right\}
< \infty. 
\end{split}
\end{equation*}

Given an initial data $ T \mapsto P_0^T $, define 
\begin{equation*}
\pi_t := P^t_0 e^{Z_t} / \E [e^{Z_t}]. 
\end{equation*}
Then, denoting the bond price in Gaussian TSM (\ref{GaussB})
by $ P^{T, \mbox{\scriptsize Gauss}}_t $, we have
\begin{equation*}
\begin{split}
& P^T_t = \E [\pi_T | \F^{p,\bW}_t ]/\pi_t \\
&= P^{T, \mbox{\scriptsize Gauss}}_t 
\exp \left\{
\int_0^t \int_E ( \delta (T,s,x) - \delta (t,s,x) ) \,\cN_p (ds, dx)
\right\} \\
&\hspace{2cm} \cdot \exp \left\{-
\int_0^t \int_E\{ ( e^{\delta (T,s,x)} -1)  -( e^{\delta (t,s,x)}-1 ) \}\ 
\nu(dx)ds
\right\}.
\end{split}
\end{equation*}
Here the filtration $ \{ \F^{p, \bW}_t \} $ is the one 
generated by $ p $ and $ \bW $.

If we further assume that 
$ \delta (t,\cdot,\cdot) $ is differentiable in $ t $ and 
if its derivative $ \partial_t \delta (t,\cdot, \cdot) $ 
is, say, uniformly bounded and if $ \nu $ is finite, 
then 
denoting the forward rate in Gaussian TSM (\ref{GaussHJM})
by $ f^{\mbox{\scriptsize Gauss}} (t,T) $, we have
\begin{equation}\label{GSHJM}
\begin{split}
&f(t,T) := - \partial_T \log P^T_t \\
&=  f^{\mbox{\scriptsize Gauss}} (t,T) 
+ \int_0^t \int_E \delta_T (T,s,x) \left\{ \cN_p (ds, dx)
-  e^{\delta (T,s,x)} \nu(dx) ds \right\}.
\end{split}
\end{equation}
This expression can be regarded as an extension of
Shirakawa's model \cite{Shi}, and 
as a special case of jump-diffusion or 
fairly general semi-martingale models,
e.g. by \cite{BjKaRu}. 

\begin{rem}
If $E$ is a finite set, 
then for each $t>0$, $\omega \in \Omega$
we define the Poisson integral of $ \delta $ as a random finite sum by 
\begin{equation*}
 \int_{E} \delta (x) N (t, dx)(\omega)
=\sum_{x \in E} \delta (x) N (t, \{ x \})(\omega),
\end{equation*}
and in this case (\ref{GSHJM}) becomes Shirakawa's model.
\end{rem}

The generalized Shirakawa model we have presented 
is basically a exponential-L\'evy model, and therefore 
we cannot use those L\'evy process without finite moments. 
For example, we miss the 
symmetric $ \alpha $-stable processes, 
which is characterized by
(constant times) $ |\xi|^\alpha $ $( \alpha \in (0,2) )$ 
as its L\'evy symbol, i.e.; the L\'evy process $ Z $ 
with the property  
\begin{equation}\label{dfn-stable}
\E [\exp \{i \langle \bxi, Z_t -Z_s \rangle ] = 
 \exp \{-(t-s) \theta |\bxi|^{\alpha}\},
\,\, \theta > 0, 
\bxi \in \R^d, 0 < s \leq t.
\end{equation}
In fact its $ \alpha $-th
moment explode. 
To construct an interest rate model driven by 
{\em stable} processes, we might use an approach 
given in the next section.

\section{L\'evy Density TSMs: 
A Generalization of QTSMs}\label{aSTSM}

\subsection{Main result}\label{aSTSM1}
Let $ Z $ be a L\'evy process in $ \R^d $
starting from the origin, adapted to a given filtration $ \{\F_t\} $.
We assume that 
it has the probability density $ p(t,x) $ with respect 
to the Lebesgue measure of $ \R^d $:
\begin{equation}\label{D1}
\P(Z_t \in dx ) = p(t,x) dx. 
\end{equation}
Here we assume that 
\begin{equation}\label{A1}
\P ( p(s, Z_t) > 0, \,\, \forall s>0, \forall t >0 ) = 1,
\end{equation}
and
\begin{equation}\label{A2}
\text{
$ p (t, \cdot ) p(s, \cdot ) $ is integrable w.r.t. 
the Lebesgue measure 
for all $ t,s > 0 $. }
\end{equation}
Then we have the following.

\begin{thm}\label{mainT}
Under (\ref{A1}) and (\ref{A2}), 
the term structure model given by
\begin{equation}\label{ALPHATSM}
P_t^T = p( \lambda_T +T-t, Z_t) /p (\lambda_t, Z_t),
\quad 0 \leq t \leq T < \infty, 
\end{equation}
where $ \lambda : [0, \infty) \to [0, \infty) $ is 
a continuous function,
is arbitrage-free in the sense of $(\ref{SPD})$ by
putting $ \pi_t = p (\lambda_t, Z_t) $.
\end{thm}

\begin{proof}
By the Markov property of $Z$, we have
\begin{equation}\label{Markovian1}
\begin{split}
\E [p (\lambda_T, Z_T)|{\cal F}_t] 
&= \E [p (\lambda_T, Z_T)|\sigma ( Z_t )] \\
&=\int_{\bf {R}^d} p( \lambda_T, x+Z_t ) \P ( Z_{T-t} \in dx) \\
&= \int_{\bf {R}^d} p( \lambda_T, x+Z_t ) p (T-t, x) \,dx.  
\end{split}
\end{equation}
Denote L\'evy symbol of $ Z $ by $ \psi $:
\begin{equation*}
\Fo [p(\lambda_T, \cdot)] (\xi) \equiv 
\int_{\R^d} e^{i \langle \xi, x \rangle} p( \lambda_T, x)  \,dx =
 e^{-\lambda_T \psi (\xi) }, \quad (\xi \in \R^d).
\end{equation*} 
Then we have 
\begin{equation*}
p (T-t, x)  = \Fo^* [e^{-(T-t) \psi (\cdot) }] (x) 
\left( \equiv (2\pi)^{-d} 
\int_{\R^d} e^{-i \langle \xi, x \rangle} e^{-(T-t) \psi (\xi) } \,d\xi \right).  
\end{equation*}
Thus we can rewrite (\ref{Markovian1}) as
\begin{equation*}
\begin{split}
\E [p (\lambda_T, Z_T)|{\cal F}_t] 
&= \langle p(\lambda_T, \cdot +Z_t), 
\Fo^* [e^{-(T-t) \psi (\cdot)}]\rangle _{L^2}. 
\end{split}
\end{equation*}
Since $ \Fo \circ \Fo^* = \mathrm{id} $ 
on the space of density functions, 
we have
\begin{equation*}
\langle p(\lambda_T, \cdot +Z_t), 
\Fo^* [e^{-(T-t) \psi (\cdot)}]\rangle _{L^2}
= (2 \pi)^{-d} \langle \Fo [p (\lambda_T,  \cdot +Z_t)], 
e^{-(T-t) \psi (\cdot)} \rangle _{L^2}. 
\end{equation*}
Observing that 
\begin{equation*}
\Fo [ p ( \lambda_T, \cdot +Z_t)] (x)
= e^{-i \langle x ,  Z_t \rangle} \Fo [p(\lambda_T, \cdot)] (x), 
\quad (x \in \R^d ),
\end{equation*}
we have 
\begin{equation}\label{Fourier1}
\begin{split}
\E [p(\lambda_T, Z_T)|{\cal F}_t]
&=  (2 \pi)^{-d}  \langle e^{- i \langle \,\cdot \, , Z_t \rangle} 
e^{-\lambda_T \psi (\cdot)}, 
{e^{- (T-t) \psi (\cdot)}} \rangle _{L^2} \\
&= (2 \pi)^{-d} \int_{\R^d}
e^{-i \langle x, Z_t \rangle } {e^{-(\lambda_T+T-t)\psi (x) }}
dx \\
&= \Fo^* [e^{-(\lambda_T+T-t)\psi (\cdot) } ] (Z_t) \\
&= p( \lambda_T+T-t, Z_t).
\end{split}
\end{equation}
Hence we have 
\begin{equation*}
\E [ p (\lambda_T, Z_T )|{\cal F}_t]/p (\lambda_t, Z_t) 
= \mbox{right-hand-side of (\ref{ALPHATSM})}.
\end{equation*}
This proves the assertion. 
\end{proof}

\subsection{Minor extensions}

\subsubsection{Translation}\label{translation}

The assumption that $ Z $ is starting from $ 0 $ is 
just a convention. 
The result is stable under the change of starting point.
More precisely, for any $ z \in \R^d $, 
substituting $ Z_t + z $ for $ Z_t $ 
in (\ref{A1}) and (\ref{ALPHATSM}) does not cause 
any problem. 

\subsubsection{Product}\label{product}
Let $ Z^1,...,Z^k $ be mutually independent 
L\'evy processes starting from the origin, $ \lambda^l :
\R_+ \to \R_+ $ be continuous maps, and 
$ p_l (t,x) $, $ l=1,\ldots,k$ be their probability density functions
satisfying (\ref{A1}) and (\ref{A2}). 
Then from Theorem \ref{mainT}, it is clear that
the bond market given by
\begin{equation*}
P^T_t = \prod_{l=1}^k p_l ( \lambda^l_T + T-t, Z^l_t )/
p_l (t, Z^l_t) 
\end{equation*}
is arbitrage-free. In fact, one can take
$ \prod_{l=1}^k p_l (t, Z_t) $ to be a state price density.

\subsection{Examples of LDTSMs}\label{exs}

\subsubsection{Quadratic Gaussian TSMs as LDTSMs}
Let $ Z $ be a $ d $-dimensional Gaussian process such that
$ Z_{t+\Delta t} - Z_t \sim N ( 0, \Sigma \Delta t )$, 
where 
$ \Sigma $ is a positive definite matrix. 
Then its probability density is given by 
\begin{equation*}
p (t,x) := ( 2 \pi  t)^{-d/2} (\det \Sigma )^{-1/2} 
\exp \{ -\frac{1}{2t}\langle \Sigma^{-1} x, x \rangle \}.
\end{equation*}
Let $ \lambda_t : \R_+ \to \R_+ $ be an increasing continuous 
function. 
Using the formula (\ref{ALPHATSM}), we have an LDTSM 
generated by $ Z $ as
\begin{equation}\label{QGas1}
\begin{split}
P^T_t &= p (\lambda_T + T-t, Z_t)/ p ( \lambda_t,Z_t) \\
& =  \left( \frac{\lambda_t}{\lambda_T + T-t} \right)^{d/2} 
\exp \{ -\frac{1}{2}\{(\lambda_T + T-t)^{-1} - \lambda_t^{-1} \}
\langle \Sigma^{-1} Z_t, Z_t \rangle \}.
\end{split}
\end{equation}
Since 
$ \bW_t := \Sigma^{-1/2} Z $ is a Brownian motion, 
(\ref{QGas1}) is also represented as
\begin{equation}\label{aff1}
P^T_t = \left( \frac{\lambda_t}{\lambda_T + T-t} \right)^{d/2} 
\exp \{ -\frac{1}{2}\{(\lambda_T + T-t)^{-1} - \lambda_t^{-1} \}
\vert \bW_t \vert^2 \}.
\end{equation} 
Note that $ X_t := \vert \bW_t \vert^2 $ is a 
$ d $-dimensional squared Bessel process, and so 
(\ref{aff1}) is a $ 1$-dimensional affine TSM. 

Set 
\begin{equation*}
A_t := \frac{1}{2\lambda_t} I \quad \text{and} \quad
k_t := -\frac{1}{2}\log ( 2 \pi  \lambda_t)^d.
\end{equation*}
Then we have 
\begin{eqnarray*}
p (\lambda_t, Z_t)
	=
\exp \left\{ 
-  \langle A_t \bW_t, \bW_t \rangle_{\R^d} + k_t  
\right\}. 
\end{eqnarray*}
Therefore, we could also obtain 
(\ref{aff1}) using the formula (\ref{QTSM1}).

It should be noted that 
a more general QTSM like (\ref{QTSM1}) 
is within our LDTSMs; it is shown 
via the "product" construction of section \ref{product}. 

\subsubsection{Cauchy TSMs}
Let $Z$ be a Cauchy process in $ \R^d $ starting from $ 0 $, 
whose probability density is given by
\begin{equation}\label{CauchyD}
p(t,x) :=\P (Z_t \in dx ) = 
\frac{ \Gamma ((d+1)/2) \theta t }{  ( \pi \{ (\theta t )^2+|x- t \gamma  |^2 \})^{(d+1)/2}},
\end{equation}
where $ \theta > 0 $ and $ \gamma \in \R^d $.
Note that Cauchy processes are strictly stable (or self-similar)
\footnote{About the stable distributions, \cite{Sat} is a good reference.}
process with parameter $ 1 $;
namely, 
\begin{equation*}
Z_{c t} \overset{\mathrm{d}}{=} c Z_t, \quad c > 0, t > 0.
\end{equation*} 
This property is seen from their L\'evy symbol: actually we have
\begin{equation*}
\E [e^{i \langle \xi, Z_t \rangle }] 
= e^{ - t ( \theta |\xi| - i \langle \xi, \gamma \rangle )}
\end{equation*}

Apparently $ p (t,x) $ in (\ref{CauchyD}) satisfies 
the assumptions  (\ref{A1}) and (\ref{A2}).
Thus we can obtain a closed form expression for an 
arbitrage-free TSM
driven by a Cauchy process as follows
\begin{eqnarray*}
P^{T}_{t}= \frac{\lambda_T + T-t}{\lambda_t}
\left( \frac{ \theta^2 \lambda_t^2+|{Z_t}- t \gamma |^2 }
{ \theta^2 (\lambda_T + T-t )^2+ | Z_t - t \gamma|^2}
\right)^{(d+1)/2}.
\end{eqnarray*}


\subsubsection{Gamma TSMs}
Let $Z$ be a one dimensional gamma process with parameters 
$a, b >0$, so that each $Z_t$ has density 
\begin{equation*}
p (t, x) =\frac{b^{at}}{\Gamma (at)} x^{at-1} e^{-bx},
\end{equation*}
for fixed $x \geq 0$. 
It is easy to check that (\ref{A1}) and (\ref{A2}) are satisfied. 
The market value at time $ t $ is given by 
\begin{eqnarray}\label{Gamma1}
\tilde{P}^{T}_{t}=
b^{\lambda_T-\lambda_t+T-t}
 \frac{\Gamma (a\lambda_t)}{\Gamma (a(\lambda_T+T-t))}
Z_t^{a(\lambda_T-\lambda_t+T-t)}.
\end{eqnarray}
Here we may think $ Z_0 > 0 $ (see section \ref{translation}). 

The expression $ \tilde{P}^T_t $ 
jumps only upwards, 
which is unrealistic for the bond price movement. 
But as we remarked in section \ref{product}, 
we can consider $ \tilde{P}^T_t $  in (\ref{Gamma1}) 
as a factor of the bond price.

\begin{rem}
We can further generalize the Gamma TSMs to include
CME (convolutions of mixtures of exponential distributions) 
and others appearing in \cite{Yam}. 
\end{rem}

\subsection*{Appendix: A derivation of (\ref{QTSM1})}
We let $ k_t \equiv 0 $ for the moment. 
By the Markov property of $\bW $, we have
\begin{equation}\label{Markovian2}
\begin{split}
& \E [\pi_T (\bW_T) |\sigma (\bW_s; s \leq t ) ] 
= \E [ \pi_T (\bW_T) |\sigma ( \bW_t )] \\
&=\int_{\bf {R}^d} \pi_T (x) \P( \bW_T \in dx | \bW_t) \\
&=\int_{\bf {R}^d} \pi_T (x+\bW_t) \P( \bW_{T-t} \in dx)\\
&=\left( \frac{1}{2\pi(T-t)} \right)^{d/2}
\int_{\bf {R}^d} \pi_T (x+\bW_t) 
\exp \left(- \frac{|x|^2}{2(T-t)}\right)dx \\
&=\left( \frac{1}{2\pi(T-t)} \right)^{d/2} 
e^{-\langle  A_T \bW_t, \bW_t \rangle} 
\int_{\bf {R}^d} e^{-2 \langle A_T \bW_t, x \rangle 
-\langle A_T x,x  \rangle_{\R^d}
- \frac{|x|^2}{2(T-t)}}dx . 
\end{split}
\end{equation}
The last equality comes from the following:
\begin{equation*}
\langle  A_T (x+\bW_t), (x+\bW_t) \rangle_{\R^d} 
= \langle  A_T \bW_t, \bW_t \rangle + 
2 \langle A_T x , \bW_t \rangle 
+\langle A_T x,x  \rangle_{\R^d}. 
\end{equation*}

Let 
$U_T$ be an orthogonal matrix such that 
$A_T = U_T^* \mathrm{diag} [\epsilon^1_T ,...,\epsilon^d_T ] U_T $.
Here, of course, $ \epsilon^1_T,..., \epsilon^d_T $ are 
eigenvalues of $ A_T $. 
Then we have
\begin{eqnarray*}
&& 2 \langle A_T x, \bW_t \rangle 
+\langle A_T x,x  \rangle_{\R^d} \\
&=& 2 \langle \mathrm{diag} [\epsilon^1_T ,...,\epsilon^d_T ] U_T x, 
U_T \bW_t \rangle_{\R^d} + 
\langle \mathrm{diag} [\epsilon^1_T ,...,\epsilon^d_T ] U_T x, 
U_T x \rangle_{\R^d}.
\end{eqnarray*}
Thus, putting $ y = (y^1_t,...,y^d_t ) = U_t x $, 
\begin{eqnarray*}
&& 2 \langle A_T x, \bW_t \rangle 
+\langle A_T x,x  \rangle_{\R^d}
+ \frac{|x|^2}{2(T-t)} \\
&=& 2 \langle \mathrm{diag} [\epsilon^1_T ,...,\epsilon^d_T ] y, 
U_T \bW_t \rangle_{\R^d} +
\langle \mathrm{diag} [\epsilon^1_T ,...,\epsilon^d_T ] y, y 
\rangle_{\R^d} +  \frac{|y|^2}{2(T-t)}, 
\end{eqnarray*}
and denoting $ c_i \equiv c_i (t,T) := \epsilon^i_T + \frac{1}{2(T-t)} $, 
\begin{eqnarray*}
&=& \sum_{i=1}^{d} \left(  c_i (y^i_t)^2 
+ 2 y^i_t (U_T \bW_t )_i \right) \\
&=&
\sum_{i=1}^{d} c_i \left\{  (y^i_t)^2 +  
2 \frac{(U_T \bW_t )_i}{c_i} + \frac{(U_T \bW_t)^2_i}{c_i^2} \right\}
- \sum_{i=1}^d \frac{(U_T \bW_t)^2_i}{c_i}.
\end{eqnarray*}
We claim here that 
\begin{equation}\label{cl1}
\sum_{i=1}^d \frac{(U_T \bW_t)^2_i}{c_i}
= 2(T-t) \langle ( 2(T-t) A_T + I )^{-1} \bW_t, \bW_t \rangle,
\end{equation}
where $ I $ is the unit matrix. 

Assuming (\ref{cl1}) for the moment, we have 
\begin{equation*}
\begin{split}
& (\text{the last expression of (\ref{Markovian2}}) ) \\
& = \left( \frac{1}{2 \pi (T-t)} \right)^{d/2}
e^{-\langle  \left\{ A_T -2(T-t) ( 2(T-t) A_T + I )^{-1}\right\} \bW_t, \bW_t \rangle} \cdot \\
& \qquad \cdot \int_{\bf {R}^d} 
e^{- \sum_{i=1}^{d} c_i \left\{  (y^i_t)^2 +  
\frac{(U_T \bW_t )_i}{c_i} \right\}^2
}
\,dy^1\cdots dy^d \\
&= \left( \frac{1}{2 (T-t)} \right)^{d/2}
e^{-\langle  \left\{ A_T - 2(T-t) ( 2(T-t) A_T + I )^{-1}\right\} \bW_t, \bW_t \rangle} (c_1 \cdots c_d)^{-1/2} \\
&= \{ \det (2 (T-t) A_T + I )\}^{-1/2}  
e^{-\langle  \left\{ A_T - 2(T-t) ( 2(T-t) A_T + I )^{-1} \right\}
\bW_t, \bW_t \rangle}. 
\end{split}
\end{equation*}
Thus we have (\ref{QTSM1}). 

Now we prove (\ref{cl1}). First observation is the following:
\begin{equation*}
\begin{split}
& \mathrm{diag} [c_1,...,c_d] 
= \mathrm{diag} [\epsilon^1_T,..., \epsilon_T^d] + \frac{1}{2(T-t)} I \\
&= U_T A_T U^*_T + \frac{1}{2(T-t)} I  
=  U_T (A_T+ \frac{1}{2(T-t)} I) U^*_T .
\end{split}
\end{equation*} 
We then notice that 
\begin{equation*}
\mathrm{diag}[\frac{1}{c_1},..., \frac{1}{c_d}]
= U_T (A_T+ \frac{1}{2(T-t)} I)^{-1} U^*_T.
\end{equation*}
Therefore, 
\begin{equation*}
\begin{split}
& \sum_{i=1}^d \frac{(U_T \bW_t)^2_i}{c_i}
= \langle \mathrm{diag}[\frac{1}{c_1},..., \frac{1}{c_d}]U_T \bW_t,
U_T \bW_t \rangle \\
& =  \langle U_T^* \mathrm{diag}[\frac{1}{c_1},..., \frac{1}{c_d}]U_T \bW_t,
\bW_t \rangle \\
& =  \langle (  A_T + \frac{1}{2(T-t)} I )^{-1} \bW_t, \bW_t \rangle \\
&= 2(T-t) \langle ( 2(T-t) A_T + I )^{-1} \bW_t, \bW_t \rangle.
\end{split}
\end{equation*}
This completes the proof.


\begin{thebibliography}{abcdefgh}

\bib{ARG}{article}{ 
author={Ahn, Dong-Hyun},
author={Dittmar, Robert F},
author={Gallant, A. Ronald},
TITLE={Quadratic Term Structure Models: Theory and Evidence},
JOURNAL={Review of Financial Studies},
date={2002},
VOLUME={15},
NUMBER={1},
PAGES={243-288},
}


\bib{AkaHar}{article}{
   author={Akahori, Jir\^o},
   author={Hara, Keisuke},
   title={Lifting quadratic term structure models 
   to infinite dimension},
   volume={16},
   number={4},
   pages={635--645},
   journal={Math. Finance},
   date={2006},
}

\bib{AlLyMa}{article}{
   author={Albeverio, Sergio},
   author={Lytvynov, Eugene},
   author={Mahnig, Andrea},
   title={A model of the term structure of interest rates based on L\'evy
   fields},
   journal={Stochastic Process. Appl.},
   volume={114},
   date={2004},
   number={2},
   pages={251--263},
   issn={0304-4149},
}


\bib{B-NO}{article}{
   author={Barndorff-Nielsen, Ole E.},
   title={Processes of normal inverse Gaussian type},
   journal={Finance Stoch.},
   volume={2},
   date={1998},
   number={1},
   pages={41--68},
   issn={0949-2984},
}

\bib{B-NL}{article}{
   author={Barndorff-Nielsen, O. E.},
   author={Levendorski{\u\i}, S. Z.},
   title={Feller processes of normal inverse Gaussian type},
   journal={Quant. Finance},
   volume={1},
   date={2001},
   number={3},
   pages={318--331},
   issn={1469-7688},
}

\bib{BjKaRu}{article}{
   author={Bj{\"o}rk, Tomas},
   author={Kabanov, Yuri},
   author={Runggaldier, Wolfgang},
   title={Bond market structure in the presence of marked point processes},
   journal={Math. Finance},
   volume={7},
   date={1997},
   number={2},
   pages={211--239},
   issn={0960-1627},
}

\bib{BMKR}{article}{
 author={Bj{\"o}rk, Tomas},
 author={Masi, G.},
   author={Kabanov, Yuri},
   author={Runggaldier, Wolfgang},
   title={Toward a general theory of bond markets},
   journal={Finance and Stochastics},
   volume={1},
   date={1997},
   pages={141--174},
}

\bib{BGM}{article}{
   author={Brace, Alan},
   author={Gatarek, Dariusz},
   author={Musiela, Marek},
   title={The market model of interest rate dynamics},
   journal={Math. Finance},
   volume={7},
   date={1997},
   number={2},
   pages={127--155},
   issn={0960-1627},
}

\bib{BriMer}{book}{
   author={Brigo, Damiano},
   author={Mercurio, Fabio},
   title={Interest rate models---theory and practice},
   series={Springer Finance},
   publisher={Springer-Verlag},
   place={Berlin},
   date={2001},
   pages={xxxviii+518},
   isbn={3-540-41772-9},
}


\bib{BroHug}{article}{
   author={Brody, Dorje C},
   author={Hughston, Lane P},
   title={Chaos and coherence: a new framework for interest-rate modelling},
   note={Stochastic analysis with applications to mathematical finance},
   journal={Proc. R. Soc. Lond. Ser. A Math. Phys. Eng. Sci.},
   volume={460},
   date={2004},
   number={2041},
   pages={85--110},
   issn={1364-5021},
}

\bib{CFP}{article}{
   author={Chen, Li},
   author={Filipovi{\'c}, Damir},
   author={Poor, H. Vincent},
   title={Quadratic term structure models for risk-free and defaultable
   rates},
   journal={Math. Finance},
   volume={14},
   date={2004},
   number={4},
   pages={515--536},
   issn={0960-1627},
}

\bib{Ch-Ni}{article}{
   author={Chiarella, Carl},
   author={Nikitopoulos Sklibosios, Christinal},
   title={A Class of Jump-Diffusion Bond Pricing Models within the HJM Framework},
   journal={Asia-Pacific Financial Markets},
   volume={10},
   date={2003},
   number={2-3},
   pages={87--127},
}

\bib{Con}{article}{
	author={Constantinides, G. M}, 
	title={A Theory of the Nominal Term Structure of Interest Rates},
	journal={Review of Financial Studies},
	volume={5},
	date={1992}, 
	pages={531--552},
}

\bib{CoTa}{book}{
   author={Cont, Rama},
   author={Tankov, Peter},
   title={Financial modelling with jump processes},
   series={Chapman \& Hall/CRC Financial Mathematics Series},
   publisher={Chapman \& Hall/CRC, Boca Raton, FL},
   date={2004},
   pages={xvi+535},
   isbn={1-5848-8413-4},
}

\bib{CIR}{article}{
   author={Cox, John C.},
   author={Ingersoll, Jonathan E., Jr.},
   author={Ross, Stephen A.},
   title={A theory of the term structure of interest rates},
   journal={Econometrica},
   volume={53},
   date={1985},
   number={2},
   pages={385--407},
   issn={0012-9682},
}

\bib{DeD}{article}{
   author={De Donno, Marzia},
   title={The term structure of interest rates as a random field: a
   stochastic integration approach},
   conference={
      title={Proceedings of 3rd Ritsumeikan Conference on Stochastic processes and applications to mathematical finance, 
},
   }, 
   book={
      publisher={World Sci. Publ., River Edge, NJ},
   },
   date={2004},
   pages={27--52},
}


\bib{DS95}{article}{
   author={Delbaen, Freddy},
   author={Schachermayer, Walter},
   title={A general version of the fundamental theorem of asset pricing},
   journal={Math. Ann.},
   volume={300},
   date={1994},
   number={3},
   pages={463--520},
   issn={0025-5831},
}

\bib{DS99}{article}{
   author={Delbaen, Freddy},
   author={Schachermayer, Walter},
   title={Non-arbitrage and the fundamental theorem of asset pricing:
   summary of main results},
   conference={
      title={Introduction to mathematical finance},
      address={San Diego, CA},
      date={1997},
   },
   book={
      series={Proc. Sympos. Appl. Math.},
      volume={57},
      publisher={Amer. Math. Soc.},
      place={Providence, RI},
   },
   date={1999},
   pages={49--58},
}

\bib{DSText}{book}{
   author={Delbaen, Freddy},
   author={Schachermayer, Walter},
   title={The mathematics of arbitrage},
   series={Springer Finance},
   publisher={Springer-Verlag},
   place={Berlin},
   date={2006},
   pages={xvi+373},
   isbn={978-3-540-21992-7},
   isbn={3-540-21992-7},
}


\bib{Duf}{book}{
   author={Duffie, Darrell},
   title={Dynamic Asset Pricing Theory},
   series={Princeton Series in Finance},
   edition={3},
   publisher={Princeton University Press},
   place={New Jersey},
   date={2001},
   pages={472},
   isbn={0-691-09022-X},
}

\bib{DFS}{article}{
   author={Duffie, Darrel},
   author={Filipovi{\'c}, Damir},
   author={Schachermayer, Walter},
   title={Affine processes and applications in finance},
   journal={Ann. Appl. Probab.},
   volume={13},
   date={2003},
   number={3},
   pages={984--1053},
   issn={1050-5164},
}

\bib{DuKa}{article}{
	author = {Duffie, Darrel},
	author={Kan, Rui},
	title = {A Yield-Factor Model of Interest Rates},
	journal = {Mathematical Finance},
	volume = {6},
	date = {1996},
	pages = {379--406},
}


\bib{EbPr}{article}{
   author={Eberlein, Ernst},
   author={Prause, Karsten},
   title={The generalized hyperbolic model: financial derivatives and risk
   measures},
   conference={
      title={Mathematical finance---Bachelier Congress, 2000 (Paris)},
   },
   book={
      series={Springer Finance},
      publisher={Springer},
      place={Berlin},
   },
   date={2002},
   pages={245--267},
}


\bib{EbRa}{article}{
   author={Eberlein, Ernst},
   author={Raible, Sebastian},
   title={Term structure models driven by general L\'evy processes},
   journal={Math. Finance},
   volume={9},
   date={1999},
   number={1},
   pages={31--53},
   issn={0960-1627},
}


\bib{EbJaRa}{article}{
   author={Eberlein, Ernst},
   author={Jacod, Jean},
   author={Raible, Sebastian},
   title={L\'evy term structure models: no-arbitrage and completeness},
   journal={Finance Stoch.},
   volume={9},
   date={2005},
   number={1},
   pages={67--88},
   issn={0949-2984},
}

\bib{Fi-Te}{article}{
   author={Filipovi{\'c}, Damir},
   author={Teichmann, Josef},
   title={On the geometry of the term structure of interest rates},
   note={Stochastic analysis with applications to mathematical finance},
   journal={Proc. R. Soc. Lond. Ser. A Math. Phys. Eng. Sci.},
   volume={460},
   date={2004},
   number={2041},
   pages={129--167},
   issn={1364-5021},
   review={\MR{2052259 (2005b:60145)}},
}


\bib{FleHug}{article}{
	author={Flesaker, B}, 
	author={Hughston, Lane P},
	title={Positive Interest},
	journal={Risk Magazine},
	volume={9},
	date={1992}, 
	number={1},
	pages={46--49},
}

\bib{HJM}{article}{
  author = {Heath, David},
  author={Jarrow, Robert},
   author={Morton, Andrew},
  title = {Bond pricing and the term structure of interest rates: a new methodology for contingent claims valuation},
  journal ={Econometrica},
  date ={1992},
  volume = {60},
  number = {1},
  pages =       {77--105}
}

\bib{HunKen}{book}{
   author={Hunt, P. J},
   author={Kennedy, J. E},
   title={Financial derivatives in theory and practice},
   series={Wiley Series in Probability and Statistics},
   edition={Revised edition},
   publisher={John Wiley \& Sons Ltd.},
   place={Chichester},
   date={2004},
   pages={xxii+437},
   isbn={0-470-86358-7},
}

\bib{HugRam}{article}{
   author={Hughston, Lane P},
   author={Rafailidis, Avraam},
   title={A chaotic approach to interest rate modelling},
   journal={Finance Stoch.},
   volume={9},
   date={2005},
   number={1},
   pages={43--65},
   issn={0949-2984},
}


\bib{IkeWat}{book}{
   author={Ikeda, Nobuyuki},
   author={Watanabe, Shinzo},
   title={Stochastic differential equations and diffusion processes},
   series={North-Holland Mathematical Library},
   volume={24},
   edition={2},
   publisher={North-Holland Publishing Co.},
   place={Amsterdam},
   date={1989},
   pages={xvi+555},
   isbn={0-444-87378-3},
}

\bib{KaWa}{article}{
   author={Kawazu, Kiyoshi},
   author={Watanabe, Shinzo},
   title={Branching processes with immigration and related limit theorems},
   journal={Teor. Verojatnost. i Primenen.},
   volume={16},
   date={1971},
   pages={34--51},
   issn={0040-361x},
}

\bib{Lev}{article}{
   author={Levendorski{\u\i}, Sergei},
   title={Pseudodiffusions and quadratic term structure models},
   journal={Math. Finance},
   volume={15},
   date={2005},
   number={3},
   pages={393--424},
   issn={0960-1627},
}

\bib{MCC}{article}{
	author={Madan, Dilip},
	author={Carr, Peter},
	author={Chang, Eric}
	title={The Variance Gamma Process and Option Pricing},
	journal={European Finance Review},
	volume={2},
	date={1998}, 
	pages={79--105},
}

\bib{Mer}{article}{
   author={Merton, Robert},
   title={Option pricing when underlying stock returns are discontinuous},
   journal={Journal of Financial Economics},
   volume={3},
   date={1976},
   pages={125--144},
}

\bib{Miy}{article}{
   author={Miyahara, Yoshio},
   title={[Geometric Levy Process \& MEMM]
   Pricing Model and Related Estimation Problems},
   journal={Asia-Pacific Financial Markets},
   volume={8},
   date={2001},
   number={1},
   pages={45--60},
}

\bib{MusRut}{book}{
   author={Musiela, Marek},
   author={Rutkowski, Marek},
   title={Martingale methods in financial modelling},
   series={Applications of Mathematics (New York)},
   volume={36},
   publisher={Springer-Verlag},
   place={Berlin},
   date={1997},
   pages={xii+512},
   isbn={3-540-61477-X},
}

\bib{Rog}{article}{
   author={Rogers, L. C. G},
   title={The potential approach to the term structure of interest rates and
   foreign exchange rates},
   journal={Math. Finance},
   volume={7},
   date={1997},
   number={2},
   pages={157--176},
   issn={0960-1627},
}

\bib{Sat}{book}{
   author={Sato, Ken-iti},
   title={L\'evy processes and infinitely divisible distributions},
   series={Cambridge Studies in Advanced Mathematics},
   volume={68},
   note={Translated from the 1990 Japanese original;
   Revised by the author},
   publisher={Cambridge University Press},
   place={Cambridge},
   date={1999},
}


\bib{Shi}{article}{
   author={Shirakawa, Hiroshi},
   title= {Interest rate option pricing with Poisson-Gaussian 
   forward rate curve processes},
   journal={Math. Finance},
   volume={1},
   date={1991},
   number={1},
   pages={77--94},
}

\bib{Vas}{article}{
	author={Vasicek, O},
	title={An equilibrium characterisation of the term structure},
	journal={Journal of Financial Economics},
	volume={5},
	date={1977},
	pages={177--188},
}

\bib{Yam}{article}{
   author={Yamazato, Makoto},
   title={Topics related to Gamma Processes},
   conference={
      title={Proceedings of 5th Ritsumeikan Conference on Stochastic processes and applications to mathematical finance},
   },
   book={
      publisher={World Sci. Publ., River Edge, NJ},
   },
   date={2006},
   pages={157--182},
}
\end{thebibliography}
\end{document}